\documentclass[reprint,amsmath,amssymb,aps,pre,showkeys,showpacs]{revtex4-1}

\usepackage{graphicx}
\usepackage{dcolumn}
\usepackage{bm}

\usepackage{color}

\newcommand{\fref}[1]{Fig.~\ref{#1}}
\newcommand{\frefs}[1]{Figs.~\ref{#1}}
\newcommand{\Fref}[1]{Figure~\ref{#1}}

\newcommand{\sref}[1]{Sec.~\ref{#1}}

\newcommand{\Sref}[1]{Section~\ref{#1}}

\newcommand{\twoeq}[2]{(\ref{#1})\&(\ref{#2})}

\def\eps{\varepsilon}

\begin{document}

%

\title{Quantitative modeling and analysis of bifurcation-induced bursting}

\author{J.E. Rubin}
\affiliation{Department of Mathematics, University of Pittsburgh, 301 Thackeray Hall, Pittsburgh, PA  15260, USA}

\author{B. Krauskopf}
\author{H.M. Osinga}
\affiliation{Department of Mathematics, University of Auckland, Private Bag 92019
Auckland 1142, New Zealand}

\date{\today}

\begin{abstract}
Modeling and parameter estimation for neuronal dynamics are often challenging because many parameters can range over orders of magnitude and are difficult to measure experimentally. Moreover, selecting a suitable model complexity requires a sufficient understanding of the model's potential use, such as highlighting essential mechanisms underlying qualitative behavior or precisely quantifying realistic dynamics. We present a novel approach that can guide model development and tuning to achieve desired qualitative and quantitative solution properties.  Our approach relies on the presence of disparate time scales and employs techniques of separating the dynamics of fast and slow variables, which are well known in the analysis of qualitative solution features. We build on these methods to show how it is also possible to obtain quantitative solution features by imposing designed dynamics for the slow variables in the form of specified two-dimensional paths in a bifurcation-parameter landscape.
\end{abstract}


\maketitle

\section{Introduction}
\label{sec:intro}
\noindent
Pulsing in lasers~\cite{Wieczorek_PhysicsReport}, mixed-mode oscillations in a chemical reaction~\cite{milik98, SNB}, and bursting in neurons~\cite{ErmentroutTerman} are just a few examples among the ubiquitous patterns that exhibit features evolving on disparate time scales. One classical approach to understanding such patterns is the method of fast-slow decomposition of the underlying mathematical model, where each variable in the model is classified as fast or slow.  All slow variables are then fixed, so that the fast variables define the so-called fast subsystem for which the slow variables are treated as parameters~\cite{pontryagin,rinzel1985, rinzel1987,mishchenko}. The slow dynamics impose a particular path that the slow variables sweep out in the bifurcation landscape of the fast subsystem. The full system dynamics can now be approximated by assuming that the fast variables drift along a sequence of attractors, punctuated by occasional rapid transitions between different attractors at certain bifurcation events, while the slow variables trace the path imposed by the slow dynamics. In particular when a single slow variable is used, the classical modeling approach considers the locus of zero speed, or nullcline of the slow variable with respect to the fast variables, and includes an additional ordinary differential equation accordingly~\cite{bertram1995, xppbook, Izhikevich}; another approach is to view the slow dynamics as a non-automomous external force~\cite{kurthsPRE2015, golubitsky2001, osingaDCDS}.

The method of fast-slow decomposition offers a simple and effective way of classifying different oscillatory patterns for two-time-scale systems. Each type of pattern is uniquely described by the sequence of bifurcations that are encountered as the slow variables trace a selected path, and by the particular attractors of the fast subsystem that are visited once such bifurcations are crossed~\cite{bertram1995, golubitsky2001, Izhikevich}. More precisely, bursting patterns will be qualitatively the same if the associated paths traced by the slow variables cross the same bifurcations of the fast subsystem in the same order. However, this approach does not reveal quantitative features of the oscillation, which may be quite important for the particular application associated with the pattern. 

In this paper, we propose a novel extension to the application of fast-slow decomposition, with two aims in mind.  First, we extend the method to the quantitative realm, so that it can be used to guide the quantitative fitting of data arising from experiments.  Second, we illustrate that fast-slow decomposition can be used as a modeling tool for systems where not enough details are known to model the dynamics of the slow variables accurately.  The underlying idea of our approach is to prescribe a simple form of slow dynamics, leading to a particular class of imposed paths for which the associated parameters can be modulated. Importantly,  we consider two slow variables, which provides access to a parameterized continuum of paths that can be used to tune quantitative features of an oscillation pattern.  The combination of imposed paths and easily obtained quantitative information about the fast subsystem can then be used to steer the fast-variable outputs to a quantitative agreement with desired targets.  The resulting system can serve as a model of the phenomenon of interest in its own right, but also as a guide for subsequent parameter estimation in the full model, if one is available, and as a tool for the development of a model for the slow components of the system, if one is not available. 

In what follows, we consider the specific field of neurophysiology and focus on patterns generated by neurons. Neurons engage in activity patterns known as bursting in a variety of settings, including sleep, novelty detection, generation of repetitive movements, release of hormones, and certain pathological conditions~\cite{burstbook,ErmentroutTerman}. In very general terms, bursting refers to any time course of the membrane potential that features active phases of consecutive high-frequency oscillations alternating with intervals in which oscillations are much smaller, much more infrequent, or absent altogether. This simple characterization, however, encompasses a striking diversity of bursting patterns that arise across different neurons in various contexts~\cite{bertram1995, Izhikevich, osingaDCDS}.  

Classic analysis of bursting dynamics, as pioneered by Rinzel~\cite{rinzel1985, rinzel1987}, utilizes the fast-slow nature of such systems and is based on singularity theory~\cite{bertram1995, golubitsky2001, osingaDCDS}. The key idea is to investigate bursting patterns and to design minimal bursting models via the analysis of the underlying fast subsystem. In particular, systems with a single slow variable have been explored extensively in this way. For such systems, an imposed path associated with a periodic bursting pattern is necessarily a line segment and the slow variable oscillates back and forth over a given range of values. The bursting patterns obtained for such systems are entirely characterized by the number and order of bifurcations of the fast subsystem that are encountered along the imposed path; for example, a complete classification of bursting patterns was attempted in~\cite{Izhikevich}.

There are cases of bursting that involve two slow variables, however. What is known as parabolic bursting is a notable case, for which the qualitative pattern of bursting has been explained in terms of slow motion across bifurcation curves of the fast subsystem~\cite{rinzelparabolic}; moreover,   a canonical underlying model associated with this bifurcation structure has been derived via nonlinear coordinate changes \cite{ermentrout1986,soto1996}. Our addition to this set of ideas is that quantitative features can be included in the design stage if two (or more) slow variables are used --- by considering an entire family of paths that encounter the required sequence of bifurcations.

As the leading example for this paper, we consider a seven-dimensional model of bursting in neurons of the respiratory brain stem, system~\cite{rubinpnas} with five fast and two slow variables. This model exhibits a particular form of neuronal bursting, which we refer to as depolarization block or DB bursting; an example of the corresponding voltage time course is shown in \fref{fig:basic}(b). This type of bursting pattern has also been observed in other neuronal models~\cite{ermentroutburst, barreto2011}. Each cycle within a DB bursting pattern consists of:\\[-5mm]
\begin{itemize}
\item[(i)] a silent phase of quiescence; \\[-7mm]
\item[(ii)] the emergence of sustained voltage spiking with gradually increasing spike frequency; \\[-7mm]
\item[(iii)] a fairly abrupt significant increase in spike frequency and attenuation of spike amplitude; \\[-7mm]
\item[(iv)] an approach towards a depolarization block state of steady, elevated voltage; \\[-7mm]
\item[(v)] a re-emergence of spiking; and \\[-7mm]
\item[(vi)] a return to quiescence.\\[-5mm]
\end{itemize}

The quantitative features of DB bursting may be important for the biological function of a neuron; for example, the spiking phase within a burst may be associated with release of substances such as hormones or with activation of a particular muscle group.  Moreover, it may not be obvious that two bursting solutions, which are qualitatively similar from the point of view of fast-slow decomposition, really share what biologists would consider to be the same features. For example, if oscillations become sufficiently small,  they will be undetectable or swamped by noise, while only certain frequencies of oscillations may suffice to achieve a biological purpose; thus, two solutions with oscillations of different amplitudes and frequencies may merit a distinction that is not present in qualitative analysis of fast-slow decomposition. 

What determines the quantitative features of a bursting pattern in a fast-slow system? We present a case study illustrating that a natural extension of the method of fast-slow decomposition can provide useful information for quantitatively fitting data on bursting arising from biological experiments.  Here, we go beyond the consideration of which bifurcation curves are crossed and take into account details of the path in the plane traced by the two slow variables in between these crossings, as well as corresponding details about the attractors of the fast subsystem encountered along this path.

This paper is organized as follows. In the next section, we present the relevant analysis of the DB bursting model from~\cite{rubinpnas} and compare its dynamics with that of a four-dimensional reduced model. The full details of the seven-dimensional model can be found in~\cite{rubinpnas}, but the system of equations is also given in the Appendix. \Sref{sec:ellipse} illustrates how imposed paths for the two slow variables in the four-dimensional model can be used to obtain quantitatively similar bursting patterns as exhibited by the full seven-dimensional model. This idea is then utilized in \sref{sec:quantify} to explore how other aspects of the fast-subsystem dynamics, which are related to but go beyond the analysis of its bifurcation diagram, determine quantitative features of the bursting patterns. We end with a discussion in \sref{sec:discussion}.

\section{The DB bursting model and its dynamics}
\label{sec:model}
\begin{figure}[t!]
  \includegraphics{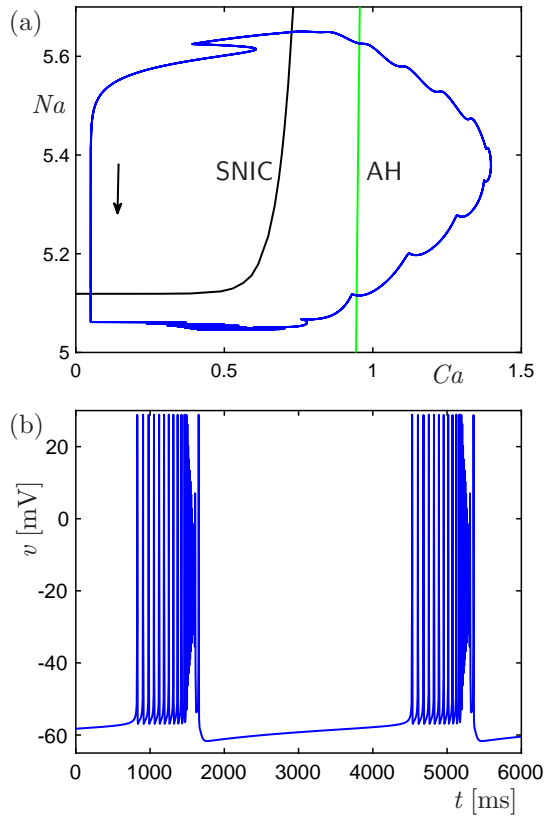}
  \caption{\label{fig:basic} 
    DB bursting in the seven-dimensional system~\twoeq{eq:7fast}{eq:7slow}. Panel~(a) shows the bursting oscillation in projection onto the $(Ca, Na)$-plane overlaid on the bifurcation set of the fast subsystem in this region, given by the loci of saddle-node bifurcations on an invariant cycle and Andronov--Hopf bifurcations, labelled {\sf SNIC} and {\sf AH}, respectively. Panel~(b) shows the corresponding time course of voltage $v$.
}
\end{figure}
\noindent
Rubin~\emph{et al.}~\cite{rubinpnas} presented a seven-dimensional model of bursting in neurons of the respiratory brain stem. The model is based on Hodgkin--Huxley formalism and is given by the following system of ordinary differential equations involving five fast and two slow variables; see the Appendix for a complete description of the various functions and parameter values used:
\begin{eqnarray}
\label{eq:7fast}
  & & \left\{ \begin{array}{rcrl}
            v' &=& -\frac{1}{c} & [I_L(v) + I_K(v, n) \\
                & & & \phantom{(} + I_{Na}(v, m, h) + I_{\rm syn}(v, s) \\
                & & & \phantom{(} + I_{\rm CAN}(v, Ca) + I_{\rm pump}(Na)], \\[1mm]
            n' &=& \frac{1}{\tau_n(v)} & [n_{\infty}(v) - n], \\[2mm]
            m' &=& \frac{1}{\tau_m(v)} & [m_{\infty}(v) - m], \\[2mm]
            h' &=& \frac{1}{\tau_h(v)} & [h_{\infty}(v) - h], \\[2mm]
            s' &=& \frac{1}{\tau_s} & [(1 - s) \, s_{\infty}(v) - k \, s], 
   \end{array} \right. \\[2mm]
\label{eq:7slow}
  & & \left\{ \begin{array}{rcrl}
            Ca' &=& \eps & [k_{\rm IP3} \, s - k_{Ca} \, (Ca - Ca_{\rm b})], \\[1mm]
            Na' &=& \alpha & [-I_{\rm CAN}(v, Ca) - I_{\rm pump}(Na)]. 
   \end{array} \right.
\end{eqnarray}
Here, system~\eqref{eq:7fast} comprises equations for the five fast variables $v$, $n$, $m$, $h$, and $s$, and system~\eqref{eq:7slow} comprises equations for the two slow variables $Ca$ and $Na$. Hence, the five-dimensional fast subsystem depends on two parameters that yield curves in the $(Ca, Na)$-parameter plane along which bifurcations occur. These curves can readily be computed with standard software packages, for example, with XPPAUT~\cite{xppbook}. 

Rubin~\emph{et al.}~\cite{rubinpnas} used fast-slow decomposition to highlight the underlying bifurcation set  of (\ref{eq:7fast}) that dictates the qualitative form of the DB bursting pattern. They found that there are two types of bifurcations, namely, a saddle-node bifurcation on an invariant cycle, denoted {\sf SNIC}, and an Andronov--Hopf bifurcation, denoted {\sf AH}; \Fref{fig:basic}(a) shows the curves  {\sf SNIC} and {\sf AH} in the $(Ca, Na)$-plane. For values $(Ca, Na)$ in between the curves {\sf SNIC} and {\sf AH}, an attracting periodic orbit of (\ref{eq:7fast}) exists togther with a single equilibrium from which it bifurcates at {\sf AH}; the equilibrium also exists to the right of {\sf AH}, where it is stable. For values $(Ca, Na)$ to the left of {\sf SNIC}, three equilibria exist; one of these is stable and corresponds to hyperpolarized voltages near the resting potential. During one period of the DB bursting pattern in the full system~\twoeq{eq:7fast}{eq:7slow}, the path taken by the slow variables $Ca$ and $Na$ is such that bifurcations are crossed four times; this is illustrated in \fref{fig:basic}(a), where the DB bursting pattern is shown in projection onto the $(Ca, Na)$-plane. Our hypothesis is that the quantitative features of this particular path taken by the two slow variables $Ca$ and $Na$ control the specific pattern of the time course of $v$ shown in panel~(b).

The stages of the DB bursting pattern listed in \sref{sec:intro} can be explained in terms of the trajectory path in the $(Ca,Na)$-plane.  Starting from a point on the periodic orbit located top-left in the $(Ca, Na)$-plane, the trajectory lies near the family of hyperpolarized rest states and the DB bursting pattern is in the quiescent or silent phase~(i); upon crossing the curve {\sf SNIC}, sustained voltage spiking emerges~(ii); as soon as the curve {\sf AH} is crossed, the frequency of the spiking increases significantly, while the amplitude decreases during stages~(iii) and~(iv); spiking re-emerges as the path crosses {\sf AH} again~(v), which is followed by a return to silence after crossing {\sf SNIC} for the second time~(vi). Hence, the active phase of each burst, namely, the epochs of voltage spiking and elevated voltage in stages (ii)--(v), encompasses the segment on the side of the curve {\sf SNIC} that contains {\sf AH}; the attenuation of spike amplitude and approach towards depolarization block of stages~(iii) and~(iv) occur between the two crossings of {\sf AH}; compare with~\cite{rubinpnas}. Barreto and Cressman~\cite{barreto2011} also associated curves of saddle-node bifurcation on an invariant cycle and Andronov--Hopf bifurcation with bursting dynamics.

For this paper, we think of the output of the seven-dimensional model as our data set, or target output for a reduced model. In an actual application, it would be replaced by experimentally recorded data.

\section{Imposed paths for a four-dimensional reduced model}
\label{sec:ellipse}
\begin{figure}[t!]
  \includegraphics{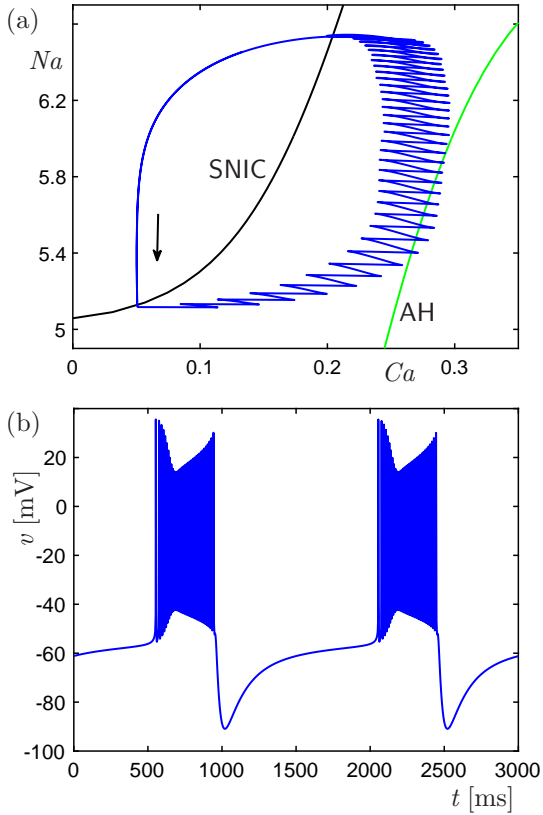}
  \caption{\label{fig:basic4D} 
    DB bursting in the reduced four-dimensional system~\twoeq{eq:4fast}{eq:7slow}. Panel~(a) shows the bursting oscillation in projection onto the $(Ca, Na)$-plane overlaid on the {\sf SNIC} and {\sf AH} curves.   Panel~(b) shows the corresponding time course of voltage $v$; compare with \fref{fig:basic}.
}
\end{figure}
\noindent
Our simulations show that a fairly similar DB-like bursting pattern can be obtained from a reduction of the seven-dimensional model of Rubin \emph{et al.}~\cite{rubinpnas} to a four-dimensional version. This reduction is achieved by a quasi-steady-state approximation, made by setting $m = m_{\infty}(v)$ and $s = s_{\infty}(v) / (s_{\infty}(v) + k)$, and by the classical step of replacing $h$ with $1 - 1.08 \, n$~\cite{krinsky, rinzel1985}. The fast subsystem is then reduced to the following two-dimensional model:
\begin{equation}
\label{eq:4fast}
  \left\{ \begin{array}{rcrl}
            v' &=& -\frac{1}{c} & [I_L(v) + I_K(v, n) \\
                & & & \phantom{(} + I_{Na}(v, m_{\infty}(v), 1 - 1.08 \, n) \\
                & & & \phantom{(} + I_{\rm syn}(v, s_{\infty}(v) / (s_{\infty}(v) + k)) \\
                & & & \phantom{(} + I_{\rm CAN}(v, Ca) + I_{\rm pump}(Na)], \\[1mm]
            n'  &=& \frac{1}{\tau_n(v)} & [n_{\infty}(v) - n].
          \end{array} \right.
\end{equation}
\Fref{fig:basic4D} shows DB bursting for the reduced system~\twoeq{eq:4fast}{eq:7slow}. The DB bursting pattern of~\twoeq{eq:4fast}{eq:7slow} is again projected onto the $(Ca, Na)$-plane in panel~(a), along with the bifurcation curves {\sf SNIC} and {\sf AH} for the two-dimensional fast subsystem~\eqref{eq:4fast}, and the time course for $v$ is shown in panel~(b). Due to the approximations used, an adjustment of some of the parameter values was required for achieving the characteristics of DB bursting of the target system~\twoeq{eq:7fast}{eq:7slow}; see the Appendix for details. One key ingredient is that we have obtained a similar bifurcation structure of the fast subsystem for~\eqref{eq:4fast} as for~\eqref{eq:7fast}.  Another is achieving the appropriate sequence of traversals of bifurcation curves.  Despite the overall similarity between the systems, comparison with \fref{fig:basic} reveals some differences with the DB bursts from the seven-dimensional model, which we found to persist despite extensive exploration of parameter space: in the DB bursts of the reduced model, stage~(ii) of sustained voltage spikes is brief, the voltage spikes are less attenuated in the approach toward depolarization block in stage~(iv), and the re-emergence of spiking is prolonged in stage~(v).  These effects are linked in \fref{fig:basic4D}(a) with the behavior of the slow variable $Ca$; the projection of the bursting solution of the four-dimensional system onto the $(Ca, Na)$-plane shows that during the spiking phase, the trajectory dances back and forth across the curve {\sf AH}, failing to cross over it completely. 

The comparison between \frefs{fig:basic} and~\ref{fig:basic4D} indicates that the model reduction from seven to four dimensions induces an unwanted effect on the behaviour of the slow variables $Ca$ and $Na$. Rather than trying to fit the slow-variable dynamics, however, our approach is to ignore their biological relevance completely. That is, we will  demonstrate that, given a particular bifurcation structure of the fast subsystem, we can generate a desired quantiative pattern by imposing a path on the associated slow variables. For this purpose, unwanted effects in a particular parameter tuning  are not relevant, as long as the bifurcation structure is present.  To this end, we replace the biological $(Ca, Na)$-dynamics with imposed paths in the $(Ca, Na)$-plane, so that the combined system has the relevant qualitative and quantitative features of DB bursting as given by the seven-dimensional model; see also~\cite{golubitsky2001, osingaDCDS} for similar ideas with one-dimensional imposed paths.

More specifically, we consider a parameterized path $\mathcal{P}$ in the $(Ca, Na)$-plane in the form of an ellipse with principal axes along the $Ca$- and $Na$-axes. Such an ellipse reasonably resembles the path shown in \fref{fig:basic}(a), yet it is defined by only five parameters that still provide enough freedom to adjust quantitative features of the dynamics. A sixth parameter is the speed $\eps$ with which the ellipse is traced. Hence, $\mathcal{P}$ is defined as
\begin{eqnarray*}
   \mathcal{P} &=& \mathcal{P}(Ca_c, Na_c, d, Ca_0, Na_0, \eps) \\
   &:=& \left\{ \begin{array}{rcr}
            Ca(t) = Ca_c &+& (Ca_0 - Ca_c) \, \cos{(\eps \, t)}\phantom{,} \\
                 &-& d \, (Na_0 - Na_c) \, \sin{(\eps \, t)}, \\[1mm]
            Na(t) = Na_c &+& (Na_0 - Na_c) \, \cos{(\eps \, t)}\phantom{.} \\
                  &+& \frac{1}{d} \, (Ca_0 - Ca_c)  \, \sin{(\eps \, t)}.
               \end{array} \right.
\end{eqnarray*}
Here, $Ca_c$ and $Na_c$ define the center of the ellipse, $d$ is its aspect ratio, and $(Ca_0, Na_0)$ is a chosen initial point at time $t = 0$. The ranges for $Ca$ and $Na$ are the intervals $[Ca_c - \delta, \; Ca_c + \delta]$ and $[Na_c - \frac{1}{d} \, \delta, \; Na_c + \frac{1}{d} \, \delta]$, respectively, where $\delta = \sqrt{(Ca_0 - Ca_c)^2 + d^2 \, (Na_0 - Na_c)^2}$. For the paths used in this paper, we fix $Na_0 = Na_c$ and use $Ca_0$ to tune the path width and $d$ to control its aspect ratio. We assume $0 < \eps \ll 1$, such that $\mathcal{P}$ is traced in the counter-clockwise direction and $Ca$ and $Na$ evolve slowly in time. 

It is convenient for our purposes that the evolution along an ellipse $\mathcal{P}$ can also be considered as the solution of the system of ordinary differential equations 
\begin{equation}
\label{eq:4slow}
  \left\{ \begin{array}{rcrl}
            Ca' &=& -\eps & d \, (Na - Na_c), \\[1mm]
            Na' &=&   \eps & \frac{1}{d} \, (Ca - Ca_c),
          \end{array} \right.
\end{equation}
for the initial conditions $Ca_0$ and $Na_0$. By combining system~\eqref{eq:4slow} with the fast subsystem~\eqref{eq:4fast} of the reduced four-dimensional model, we obtain the four-dimensional \emph{driven} model ~\twoeq{eq:4fast}{eq:4slow} that is the main subject of this paper. Note that the formulation~\eqref{eq:4slow} suffices for our purposes here, even though its periodic orbits --- the ellipses --- are not isolated. In case continuation of the overall periodic orbit of the driven system ~\twoeq{eq:4fast}{eq:4slow} is required, a family of isolated attracting elliptical periodic orbits can be obtained in similar fashion from a scaled version of the Hopf normal form~\cite{xppbook}.

\section{Results}
\label{sec:quantify}
\noindent
The approach of imposing a path $\mathcal{P}$ for the two slow variables $Ca$ and $Na$ allows us to control and  explore their qualitative and quantitative effects systematically. Our main result is that we are able to quantify different features of the DB bursting patterns in terms of the underlying fast subsystem; these findings go beyond the standard approach of determining the bifurcation diagram of the fast subsystem.  

\subsection{Bursting with an imposed $(Ca, Na)$-ellipse}
\begin{figure}[t!]
  \centering
  \includegraphics{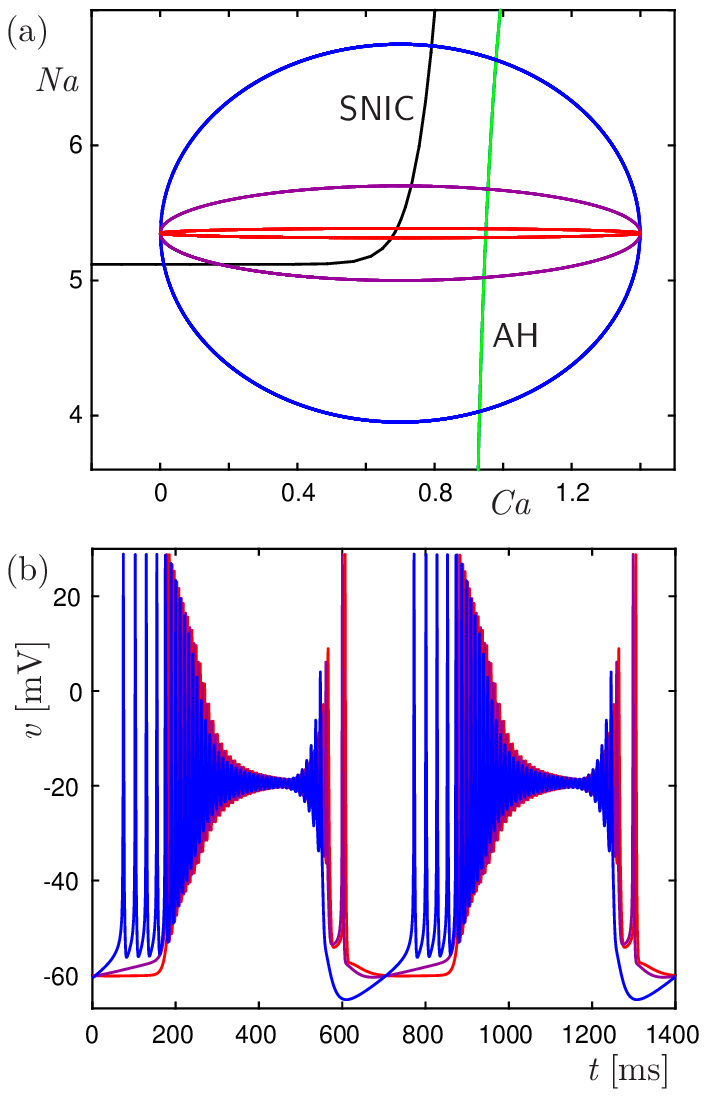}
  \caption{\label{fig:burst} 
    DB bursting patterns for the seven-dimensional driven system~\twoeq{eq:7fast}{eq:4slow} with $Ca_c = 0.7$, $Na_c = 5.35$, $Ca_0 = 0$ and $\eps = 0.009$, and three different values for $d$, namely, $d = \frac{1}{2}$ (blue), $d = 2$ (purple) and $d = 20$ (red). Panel~(a) shows the three paths $\mathcal{P}$ in the $(Ca, Na)$-plane overlaid on the loci for {\sf SNIC} and {\sf AH}, and panel~(b) shows the three corresponding time courses for $v$.}
\end{figure}
\noindent
We expect that a driven system, where the two slow variables $Ca$ and $Na$ evolve as determined by an imposed elliptic path, can exhibit DB bursting if the fast subsystem gives rise to a region in the $(Ca, Na)$-plane bounded by bifurcation curves {\sf SNIC} and {\sf AH}. We now confirm that claim by showing that the qualitative features of DB bursting are obtained by choosing a path $\mathcal{P}$ that crosses these curves in the specific sequence {\sf SNIC}, {\sf AH}, {\sf AH}, and {\sf SNIC}.

We first explore this for a driven version of the seven-dimensional model, that is, by combining the fast subsystem~\eqref{eq:7fast} with the driven system~\eqref{eq:4slow}. \Fref{fig:burst} shows three different DB bursting patterns for system~\twoeq{eq:7fast}{eq:4slow}, using three different imposed paths for the slow variables $Ca$ and $Na$. The imposed paths in the $(Ca, Na)$-plane are shown in panel~(a) overlaid on the loci for {\sf SNIC} and {\sf AH} of the fast subsystem~\eqref{eq:7fast}; the corresponding three time courses for $v$ are shown in panel~(b). For each path, we used $Ca_c = 0.7$, $Na_c = 5.35$, with initial condition $Ca_0 = 0$ (recall that we always set $Na_0 = Na_c$) and constant speed $\eps = 0.009$. Only the aspect ratio $d$ was varied, namely, $d = \frac{1}{2}$ (blue) for the fattest ellipse shown in panel~(a), which was increased to $d = 2$ (purple) and then $d = 20$ (red). All three paths cross the curves {\sf SNIC} and {\sf AH} in the specified sequence needed for DB bursting and, as illustrated in panel~(b), all yield DB bursting patterns that exhibit the qualitative features described for stages~(i)--(vi); note that the path with $d = \frac{1}{2}$ yields some tonic spiking at the onset of each active phase that is also seen in the DB bursting of system~\twoeq{eq:7fast}{eq:7slow}, and which is a quantitative feature not seen for the other paths.

\begin{figure}[t!]
  \centering
  \includegraphics{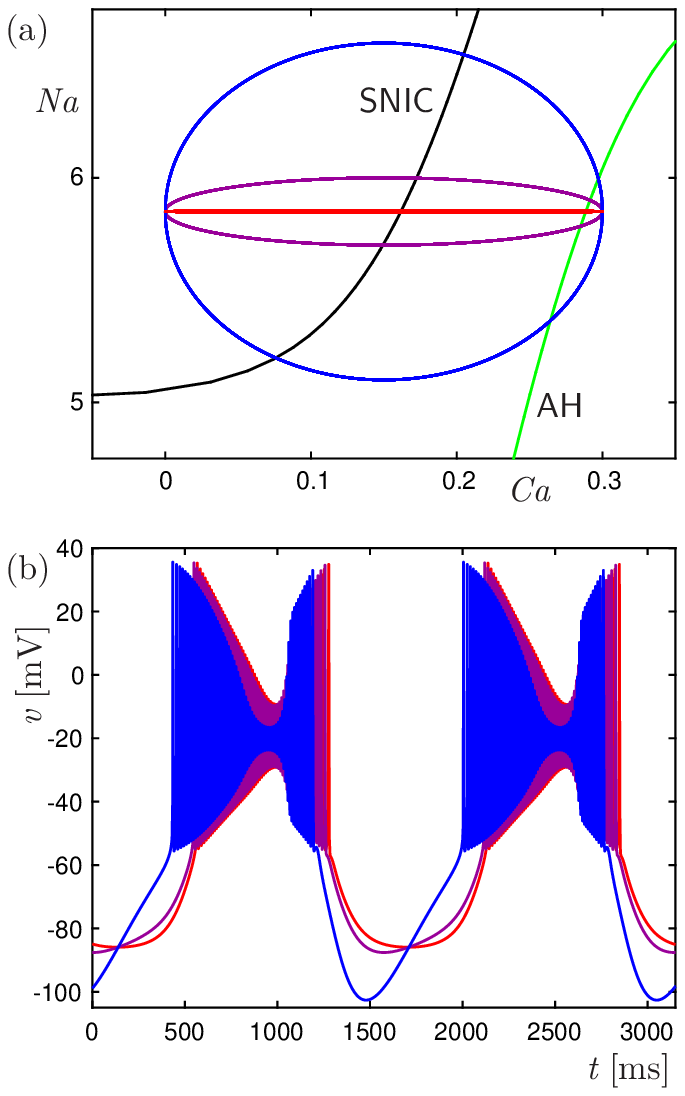}
  \caption{\label{fig:burst4D} 
    DB bursting patterns for the four-dimensional driven system~\twoeq{eq:4fast}{eq:4slow} with $Ca_c = 0.15$, $Na_c = 5.85$, $Ca_0 = 0$ and $\eps = 0.004$, and three different values for $d$, namely, $d = \frac{1}{5}$ (blue), $d = 1$ (purple) and $d = 50$ (red). Panel~(a) shows the three paths $\mathcal{P}$ in the $(Ca, Na)$-plane overlaid on the loci for {\sf SNIC} and {\sf AH}, and panel~(b) shows the three corresponding time courses for $v$.}
\end{figure}
Qualitatively similar DB bursting can be obtained for the four-dimensional driven system~\twoeq{eq:4fast}{eq:4slow}. Since the loci of {\sf SNIC} and {\sf AH} for the fast subsystem~\eqref{eq:4fast} have shifted slightly in the $(Ca, Na)$-plane, slightly different paths must be imposed. \Fref{fig:burst4D} shows the effect of three different ellipses centered at $Ca_c = 0.15$ and $Na_c = 5.85$, with $Ca_0 = 0$ and $\eps = 0.004$. Again, only the aspect ratio $d$ was varied, where we chose $d = \frac{1}{5}$ (blue), $d = 1$ (purple) and $d = 50$ (red). As for \fref{fig:burst}, the associated paths are overlaid in panel~(a) on the loci for {\sf SNIC} and {\sf AH} in the $(Ca, Na)$-plane, and the three corresponding time courses for $v$ are shown in panel~(b). Note that only a short segment of each of the three imposed elliptic paths lies in the region to the right of {\sf AH}, which regulates stages~(iv) and~(v) of DB bursting. Consequently, the approach towards and away from the depolarization block state of steady, elevated voltage is relatively brief; compare \frefs{fig:burst}(b) and \ref{fig:burst4D}(b). Also, the path with $d = \frac{1}{5}$ (blue) yields particularly deep hyperpolarizations between bursts, and maximal spike attenuation during bursts, which is also the case for the path with $d = \frac{1}{2}$ (blue) in \fref{fig:burst}. The difference in the range for the $t$-axes in \frefs{fig:burst}(b) and \ref{fig:burst4D}(b) is entirely controlled by the different choices for $\eps$.

Thus, as expected, the qualitative features of the DB bursting pattern for both the full seven-dimensional driven system~\twoeq{eq:7fast}{eq:4slow} and the reduced four-dimensional driven system~\twoeq{eq:4fast}{eq:4slow} are entirely determined by the type of bifurcations exhibited by their respective fast subsystems and the specific order in which these bifuractions are encountered. As illustrated in \frefs{fig:burst} and~\ref{fig:burst4D}, we can easily control such properties by choosing $Na_c$ large enough and deciding on appropriate values for $Ca_c$ and $Ca_0$. We also observe, however, that quantitative differences in the bursting patterns are obtained by tuning the aspect ratio $d$ of the elliptic path, as well as the speed $\eps$ at which this path is traced.  For example, comparison between \frefs{fig:burst} and~\ref{fig:burst4D} suggests that the speed with which $Ca$ and $Na$ traverse the same path leads to quantitative differences in the bursting.  We next pursue the idea of linking model parameters associated with the imposed path in $(Ca,Na)$ space to quantitative features of the DB pattern.

\subsection{Quantitative features of bursting beyond the basic bifurcation structure}
\noindent

We now explore how the dynamics of the fast subsystem can be analysed to determine quantitative features of bursting solutions for the full system. Here, we focus exclusively on the reduced four-dimensional driven system~\twoeq{eq:4fast}{eq:4slow}.

\begin{figure}[t!]
  \centering
  \includegraphics{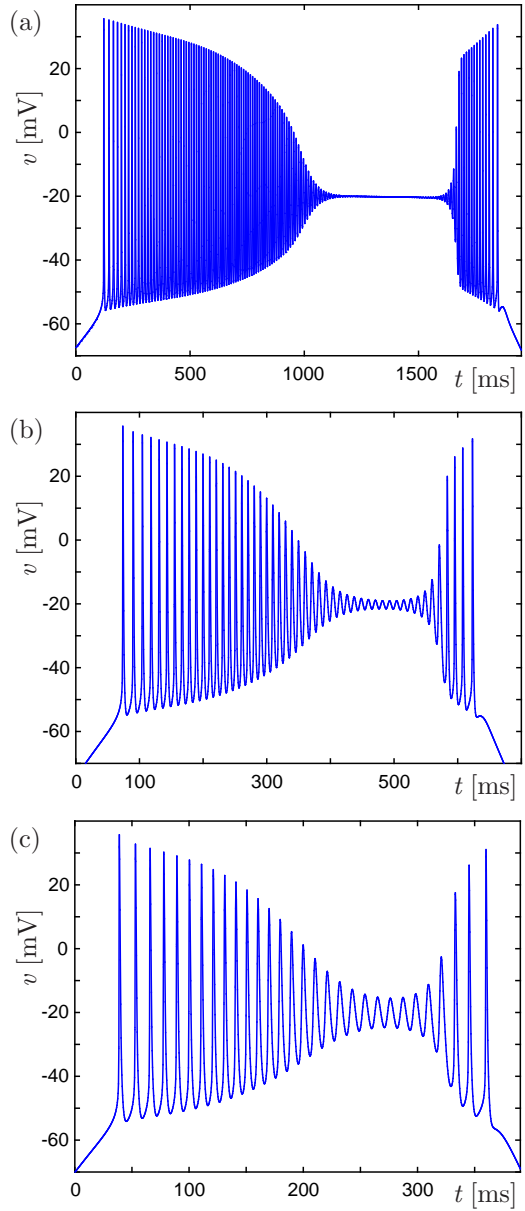}
  \caption{\label{fig:vareps} 
    The same path in the $(Ca, Na)$-plane can yield quantitative differences in bursting features, depending on the speed at which it is traversed.  Shown are the time courses of $v$ for the four-dimensional driven system~\twoeq{eq:4fast}{eq:4slow}  resulting from the same imposed elliptic path with speeds $\eps = 0.002$ in panel~(a), $\eps = 0.006$ in panel~(b), and $\eps = 0.01$ in panel~(c). The other parameters for the path are $Ca_c = 0.15$, $Na_c = 5.85$, $Ca_0 = 0$, and $d = 0.1$.}
\end{figure}

We first return to the consideration of the traversal rate of the imposed path.  \Fref{fig:vareps} shows time courses of $v$ for three quantitatively different bursting patterns for the four-dimensional driven system~\twoeq{eq:4fast}{eq:4slow} with $\eps = 0.002$ in panel~(a), $\eps = 0.006$ in panel~(b), and $\eps = 0.01$ in panel~(c) for an imposed elliptic path with the parameters $Ca_c = 0.15$, $Na_c = 5.85$, $Ca_0 = 0$ and $d = 0.1$. Consistent with the previous subsection, when $\eps$ is small, each phase of the burst is elongated due to the smaller traversal rate; note the different $t$-scales in the three panels. Moreover, as $\eps$ increases, the faster traversal of the same path in the $(Ca, Na)$-plane yields less attenuation of the voltage oscillations during the depolarization block phase of each burst cycle. As is particularly clear in \fref{fig:vareps}(a), a slower traversal intensifies this attenuation.

We dig deeper into the source of quantitative differences in patterns by next considering the time intervals between spikes in stage~(ii) of a DB bursting pattern, during which repetitive, high-amplitude voltage spikes occur. These spikes arise from a family $\Gamma_0$ of stable periodic orbits of the fast subsystem that exist in the $(Ca, Na)$-region in between the bifurcation curves {\sf SNIC} and {\sf AH}. The time intervals between spikes can vary significantly across successive spike pairs within a bursting solution, and the time until emergence of the first spike after th esecond crossing of {\sf SNIC} can differ appreciably across parameter values and paths. Our hypothesis is that these quantitative features observed in stage~(ii) of DB bursting are controlled by characteristics of the $(Ca, Na)$-dependent family $\Gamma_0$.   

\begin{figure}[t!]
\begin{center}
\includegraphics{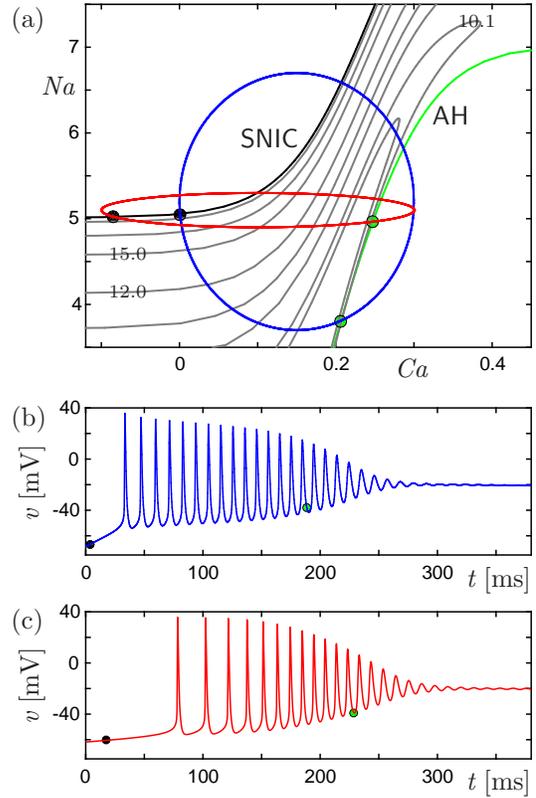}
\caption{\label{fig:period} 
  Interspike intervals of the DB bursting pattern are significantly affected by the different periods of the periodic orbits in the family $\Gamma_0$ of the fast subsystem that are encountered along the imposed path. Panel (a) shows {\sf SNIC} and {\sf AH} in the $(Ca, Na)$-plane along with contours (grey) of equal period for periodic orbits in $\Gamma_0$. Overlaid are a vertically wide (blue) path with $Ca_c = 0.15$, $Na_c = 5.2$, $d = 0.1$, and $Ca_0 = 0$ and a horizontally elongated (red) path with $Ca_c = 0.1$, $Na_c = 5.1$, $d = 1$, and $Ca_0 = -0.1$. The corresponding time courses of $v$ over the first half period are shown in panels~(b) and~(c), respectively; their crossings of the curves {\sf SNIC} and {\sf AH} are marked with green and black circles, respectively.
}
\end{center}
\end{figure}
Specifically, we note that each periodic orbit in $\Gamma_0$ has a well-defined period that depends on the choice for $Ca$ and $Na$. The periods encountered along a path $\mathcal{P}$ are determined by the location of this imposed path in the region bounded by {\sf SNIC} and {\sf AH}. We find that the progression of periods encountered along $\mathcal{P}$ strongly shapes the time intervals between spikes in stage~(ii) of the resulting DB bursting pattern. \Fref{fig:period}(a) shows the curves {\sf SNIC} and {\sf AH} in the $(Ca, Na$)-plane together with contours of equal period for the periodic orbits of the family $\Gamma_0$.
The periods associated with the contours progressively decrease as we move farther from {\sf SNIC},  starting from 40 for the curve closest to the SNIC, although close to {\sf AH} some non-monotonicity sets in (e.g, note the turning point for the labeled contour with period 10.1). This panel also includes two elliptic paths with $\eps = 0.009$: the fat (blue) path is given by $Ca_c = 0.15$, $Na_c = 5.2$, $d = 0.1$, and $Ca_0 = 0$, and the thin (red) path by $Ca_c = 0.1$, $Na_c = 5.1$, $d = 1$, and $Ca_0 = -0.1$. The corresponding time traces of $v$ over the first half period are shown in panels~(b) and~(c) of \fref{fig:period}, respectively. The filled black and green circles indicate the moment when the solutions cross {\sf SNIC} and {\sf AH}, respectively. 

The fat (blue) path crosses {\sf SNIC} roughly orthogonal to the period contours and thus progresses quickly to a region with relatively low periods.  The interspike intervals in the corresponding time trace of $v$, shown in \fref{fig:period}(b), are relatively short, with little change over the period of the burst. In contrast, the thin (red) path crosses {\sf SNIC} in a direction that is aligned with the high-period contours close to the curve {\sf SNIC}.  Correspondingly, there is a longer delay from the crossing (black circle) to the first spike in the time trace of $v$, shown in \fref{fig:period}(c), and the first few interspike intervals are longer than any of the ones shown in \fref{fig:period}(b), with a gradual compression of spike times over the course of the burst.  Note that we selected parameter values such that the overall times from the start (black circles) to the end (green circles) of the active phase, in between {\sf SNIC} and {\sf AH}, are quite similar for both DB bursts, so the longer interspike intervals in \fref{fig:period}(c) are not simply due to a slower passage through stage~(ii) of the DB bursting pattern. 

\begin{figure}[t!]
\begin{center}
\includegraphics{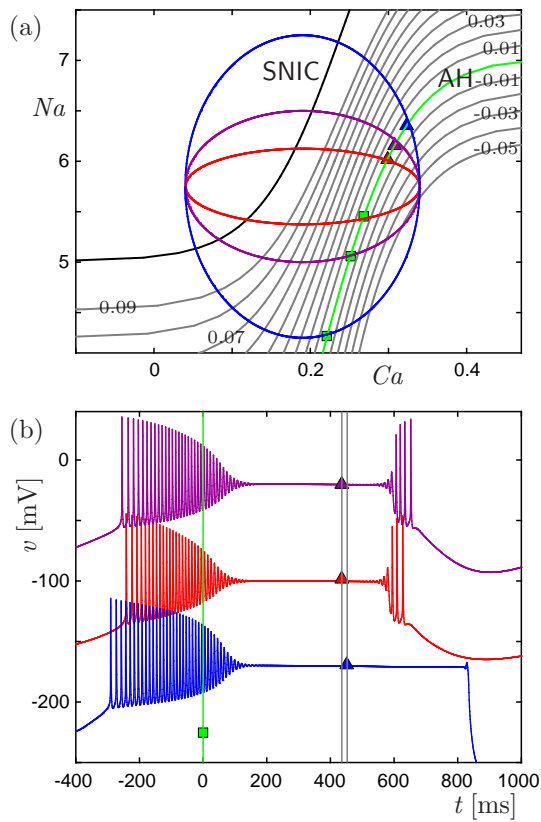}
\caption{\label{fig:nonmon} 
  Contraction to and expansion away from the depolarization block state depends on the real part $Re(\lambda)$ of the eigenvalues of the associated equilibria in the fast subsystem that are encountered along the imposed path. Panel (a) shows {\sf SNIC} and {\sf AH} in the $(Ca, Na)$-plane along with contours (grey) of $Re(\lambda)$. Overlayed are paths with $Ca_c = 0.19$, $Na_c = 5.75$, $Ca_0 = 0.04$ and $\eps = 0.004$, and varying $d = 0.1$ (blue), $d = 0.2$ (purple) and $d = 0.4$ (red). The corresponding time courses of $v$ are shown in panel~(b) with time shifted so that the first crossing occurs at $t = 0$ (green square), and voltage shifted down by $200$ and $100 \, {\rm mV}$ for the time courses with $d = 0.1$ (blue) and $d = 0.4$ (red), respectively; the second crossing of {\sf AH} is marked with triangles. 
}
\end{center}
\end{figure}
Next, we turn to stages~(iv) and~(v) of the DB bursting pattern and quantify the approaches towards and away from the depolarization block state of steady, elevated voltage, which occur after the first and second crossing of {\sf AH}, respectively.  The observed rate of contraction (expansion) is related to the stable (unstable) eigenvalues obtained from the linearization of the fast subsystem about the equilibria associated with this elevated state. Near the curve {\sf AH} of Andronov--Hopf bifurcations, these equilibria all have complex-conjugate pairs of eigenvalues. Their real parts are negative to the right of {\sf AH} and positive to the left of it. \Fref{fig:nonmon}(a) shows $15$ contours in the $(Ca, Na)$-plane of these real parts $Re(\lambda)$, uniformly distributed from $-0.05$ to $0.09$; the curve {\sf AH} corresponds to the countour $Re(\lambda) = 0$. We select three elliptic paths that are all centered at $Ca_c = 0.19$ and $Na_c = 5.75$ and use $Ca_0 = 0.04$ and $\eps = 0.004$. The only difference in these contours is in their aspect ratios, which are $d = 0.1$ (blue), $d = 0.2$ (purple) and $d = 0.4$ (red). Their intersections with {\sf AH} are marked with squares and triangles, indicating the first and second crossings, respectively. The corresponding time courses for $v$ are shown in \fref{fig:nonmon}(b), where time is shifted so that the first crossing occurs at $t = 0$ (green square).   Here, the burst with $d = 0.2$ (purple, top) is plotted in actual coordinates, while the ones with $d = 0.1$ (blue, bottom) and $d = 0.4$ (red, middle) have been shifted down by $200$ and $100 \, {\rm mV}$, respectively. 

As shown in \fref{fig:nonmon}(b), stage~(iv) of the DB bursting patterns is rather similar across these parameter values, even though the times until the second crossing of {\sf AH} (triangles) are different. Observe that the times between AH crossings for the bursts with $d = 0.2$ (purple, top) and $d = 0.4$ (red, middle) are similar, while the burst with $d = 0.1$ (blue, bottom) exhibits a significantly longer stage~(iv). \fref{fig:nonmon}(a) illustrates that the burst with $d = 0.4$ (red) does not cross the contour $Re(\lambda) = -0.05$, which indicates that the net contraction towards the depolarization block state is not as strong as for the other two bursts. Consequently, even though the bursts with $d = 0.2$ (purple) and $d = 0.4$  (red) spend about the same time to the right of {\sf AH}, the time needed to expand away from the depolarization block state is shorter for the burst with $d = 0.4$  (red), resulting in a shorter time until the onset of stage~(v); see \fref{fig:nonmon}(b).  Note further that the longer time spent to the right of {\sf AH} for the burst with $d = 0.1$ (blue) results in a stronger contraction towards the depolarization block state and subsequent slower expansion away from it. Hence, this case includes virtually no additional spikes before  the initiation of the silent phase at the end of stage~(v). 

In fact, the times spent between crossings of bifurcation curves and the rates of contraction and expansion encountered during these times, as indicated by the contours in \fref{fig:nonmon}(a), combine in such a way that the total number of additional spikes exhibited in stage~(v), before the return to the silent phase, depends non-monotonically on the aspect ratio of the imposed elliptic path.  The time courses in \fref{fig:nonmon}(b) are ordered to emphasize this non-monotonicity; they reveal that the largest number of additional spikes occur for the (purple) elliptic path with $d = 0.2$.  Although these spikes start earlier within the burst, for larger $d$, such as $d = 0.4$ (red), the relatively short time from the second {\sf AH} crossing to the second {\sf SNIC} crossing limits the number of spikes that can occur within stage~(v).  For smaller $d$, such as $d=0.1$ (blue), the need for additional expansion to overcome the strong contraction from stage~(iv) is so strong that not even a full spike can be fired within stage~(v) before {\sf SNIC} is crossed and the active phase ends.

\section{Discussion}
\label{sec:discussion}
\noindent
We introduced and demonstrated a conceptual approach to modeling and analysis of fast-slow dynamics that allows one to determine and extract quantitative information of relevance to the appplication at hand. Our approach applies to fast-slow systems with at least two slow variables.  The qualitative features of the system dynamics are then given by the bifurcation structure of the fast subsystem and how it is encountered by a closed path in the space of slow variables. Obtaining information on the fast subsystem in the regions that are traversed along this path, such as periods of periodic orbits and contraction rates associated with invariant objects, allows one to characterize how specific quantitative features of the observed solution arise. Moreover, we specify a family of paths in the state space of slow variables suitable for the construction of an associated driven system, which we use to investigate systematically these quantitative solution characteristics.

As our specific test-case example we considered depolarization block or DB bursting as displayed by a seven-dimensional model  from~\cite{rubinpnas} for neurons of the respiratory brain stem. Here, the concentrations $Ca$ of calcium and $Na$ of sodium in the neuron evolve on a much slower time scale than the voltage potential $v$ across the cell membrane, which is the observable in experiments. We derived a four-dimensional reduction of that model with only two fast variables and also considered its driven version, where the family of paths was chosen to be ellipses in the plane of $Ca$ and $Na$.  DB bursting is organized qualitatively by a curve of saddle-node bifurcations on an invariant cycle ({\sf SNIC}) and a curve of Andronov--Hopf bifurcations ({\sf AH}) of the fast subsystem, and elliptic paths naturally capture the crossing sequence arising in the biological model. The bursting pattern traces a family of attracting periodic orbits of the fast subsystem, and also exhibits contraction towards and subsequent expansion away from the equibria that undergo the Andronov--Hopf bifurcation.  Although there are some quantitative discrepancies between the active phase dynamics in the seven- and four-dimensional biological models, a key point is that by tuning imposed paths, we can achieve a more complete quantitative agreement.  

More specifically, we found that the amount of time an imposed path spends in different regions of the $(Ca, Na)$-plane, corresponding to passage between bifurcation curve crossings, is only one factor in determining timing-related quantitative features of the resulting bursting pattern. In fact, the periods of the periodic orbits of the fast subsystem that are encountered along the path play a crucial role in shaping the spike frequency within the bursts; this relation can be analyzed by computing contours of equal period in the region where (attracting) periodic orbits of the fast subsystem exist. Moreover, we found that abrupt transitions from slow spikes (at the onset of spiking) to faster spikes can be predicted this way. Similarly, the eigenvalues of the equilibria of the fast subsystem that are encountered along the path are responsible for the observed contraction and expansion rates associated with the depolarization block state of DB bursting; details can be analyzed by computing contours of equal real parts of such eigenvalues in the relevant regions of the $(Ca, Na)$-plane. These contours explain an apparent and initially counter-intuitive non-monotonicity in the number of spikes that arise at the end of the DB phase, just before the onset of the silent phase, with respect to the aspect ratio $d$ of the imposed elliptical path in the $(Ca, Na)$-plane.  We find similar results (not reported here) with the models proposed in~\cite{barreto2011} and in~\cite{ermentroutburst}. Given their connection to the well-established relationship between bifurcation structures of the fast subsystem and qualitative features of the bursting pattern, we believe that our ideas will extend naturally to models with different underlying bifurcation diagrams of their associated fast subsystems; for example, our approach might be of use for modeling and analysis of the sound patterns produced by songbirds~\cite{amador2014, amador2015}.  

The linkage that we have established between quantitative bursting features and path of traversal through the $(Ca, Na)$-plane opens the door to designing an imposed path for the slow variables $Ca$ and $Na$ to achieve a particular quantitative outcome. This ability to control solution features is very useful from a model development perspective, especially in situations where it is advantageous to work with a lower-dimensional model reduction that must be tuned to capture complicated dynamics. The design of paths can guide model development and, moreover, it can be used for parameter estimation and fitting for the slow dynamics,  which may be particularly useful in light of difficulties in experimentally measuring quantities associated with slow variables in neural models. Once a family of paths has been chosen, the path or paths with the most suitable quantitative features can be selected. We used ellipses here, but other, more complicated paths may be considered, for example, to achieve a closer fit between the voltage dynamics of the reduced four-dimensional model and the target seven-dimensional model. An important point is that the family of paths is specified by a resonably low number of parameters. In future work, we plan to consider the utility of our proposed approach for parameter estimation of fast-slow bursting models by coupling it with optimization techniques. One challenge is to define suitable functions that encapsulate, in terms of relevant quantitative features, the distance of a given (periodic) output from the dynamics under consideration; in applications, the latter may be generated by an actual experiment, rather than a higher-dimensional mathematical model.

\begin{acknowledgments}
This work was partially supported by NSF awards DMS 1312508 and DMS 1612913.  JR thanks the University of Auckland for hosting an extended visit while part of this work was performed.  
\end{acknowledgments}

\appendix
\section*{Appendix}
\noindent
For completeness, we provide here the precise definitions of the functions used in system~\twoeq{eq:7fast}{eq:7slow}; the parameters are given in Table~\ref{tab:7dpar}. Recall that system~\twoeq{eq:7fast}{eq:7slow} is the same model of DB bursting in neurons of the respiratory brain stem as described in~\cite{rubinpnas} and full details can also be found there. 

The equation for voltage $v$ in system~\eqref{eq:7fast} includes the following currents:
\begin{displaymath}
  \begin{array}{rcl}
    I_L(v)                 &=& g_L \, (v - E_L), \\[1mm]
    I_K(v, n)             &=& g_K \, n^4 \, (v - E_K), \\[1mm]
    I_{Na}(v, m, h)     &=& g_{Na} \, m^3 \, h \, (v - E_{Na}), \\[1mm]
    I_{\rm syn}(v, s)    &=& g_{\rm syn} \, s (v - E_{\rm syn}), \\[1mm]
    I_{\rm CAN}(v, ca) &=& {\displaystyle \frac{g_{\rm CAN} \, (v - E_{\rm CAN})}{1 + \exp((Ca-k_{\rm CAN})/\sigma_{\rm CAN})}}, \\[3mm]
    I_{\rm pump}(Na)   &=& r_{pump} \, (\phi(Na) - \phi(Na_b)), 
  \end{array}
\end{displaymath}
where
\begin{displaymath}
  \phi(Na) = \frac{Na^3}{Na^3 + k_{Na}^3}. 
\end{displaymath}
The equations for $n$, $m$, $h$, and $s$ all have a similar form, with 
\begin{displaymath}
  X_{\infty}(v) = \frac{1.0}{1.0 + \exp((v-\theta_X)/\sigma_X)},
\end{displaymath}
where $X \in \{n, m, h, s\}$, and 
\begin{displaymath}
  \tau_X(v) = \frac{t_X}{\cosh\left(\frac{v - \theta_X}{2 \, \sigma_X}\right)},
\end{displaymath}
where $X \in \{n, m, h\}$; the time scale $\tau_s(v) = \tau_s$ for $s$ is taken constant. There are many parameters in this model, and their values as used in this paper are listed in Table~\ref{tab:7dpar}.
\begin{table}[h!]
  \centering
  \caption{\label{tab:7dpar}
    Parameter values used for the functions in the seven-dimensional model~\twoeq{eq:7fast}{eq:7slow}.}
\begin{tabular}{|r@{ $=$ }r@{.}l|r@{ $=$ }r@{.}l|r@{ $=$ }r@{.}l|} 
\hline 
\hline 
\multicolumn{3}{|c|}{conductances} & \multicolumn{3}{c|}{reversal} & \multicolumn{3}{c|}{half} \\
\multicolumn{3}{|c|}{(${\rm nS}$)} & \multicolumn{3}{c|}{potentials (${\rm mV}$)} & \multicolumn{3}{c|}{activations} \\
\hline
$g_L$    &     $3$ & $0$ & $E_L$     & $-60$ & $0$      & $\theta_h$ & $-30$ & $0 \, {\rm mV}$ \\
$g_{Na}$ & $150$ & $0$ & $E_{Na}$ &   $85$ & $0$      & $\theta_m$ & $-36$ & $0 \, {\rm mV}$ \\
$g_{K}$  &   $30$ & $0$ & $E_{K}$   & $-75$ & $0$      & $\theta_n$ & $-30$ & $0 \, {\rm mV}$ \\ 
$g_{\rm syn}$ & $2$ & $5$ & $E_{\rm syn}$ & $0$ & $0$  & $\theta_s$ & $15$ & $0 \, {\rm mV}$ \\
$g_{\rm CAN}$ & $4$ & $0$ & $E_{\rm CAN}$ & $0$ & $0$ & $k_{\rm CAN}$ & $0$ & $9 \, \mu{\rm M}$ \\
\hline 
\hline 
\multicolumn{3}{|c|}{slopes} & \multicolumn{3}{c|}{time} & \multicolumn{3}{c|}{} \\
\multicolumn{3}{|c|}{(${\rm mV}$ or $\mu{\rm M}$)} & \multicolumn{3}{c|}{constants (${\rm ms}$)} & \multicolumn{3}{c|}{scaling constants} \\
\hline
$\sigma_h$           & $5$         & $0$ & $t_h$       & $15$ & $0$ & $k_{Na}$ & $10$ & $0\phantom{0} \, {\rm mM}$ \\
$\sigma_m$           & $-8$      & $5$ & $t_m$       &   $1$ & $0$ & $Na_b$ & $5$ & $0\phantom{0} \, {\rm mM}$ \\
$\sigma_n$           & $-5$       & $0$ & $t_n$        & $30$ & $0$ & $k_{Ca}$ & $22$ & $5\phantom{0} \, {\rm ms}^{-1}$ \\
$\sigma_s$           & $-3$        & $0$ & $\tau_s$ & $15$ & $0$ & $Ca_b$ & $0$ & $05 \, \mu{\rm M}$ \\
$\sigma_{\rm CAN}$ & $-0$ & $05$ & \multicolumn{3}{c|}{}   & $k_{\rm IP3}$ & \multicolumn{1}{l}{$1200$} & \multicolumn{1}{r|}{$\mu{\rm M ms}^{-1}$} \\
\hline
\hline
\multicolumn{3}{|c}{other} & \multicolumn{6}{c|}{} \\
\hline
$C$ & $45$ & $0 \, {\rm pF}$ & $k$ & \multicolumn{2}{l|}{$1.0$} & $r_{\rm pump}$ & $200$ & $0 \, {\rm pA}$ \\
\hline
$\eps$ & \multicolumn{2}{l|}{$7 \times 10^{-4}$} & $\alpha$ & \multicolumn{5}{l|}{$6.6 \times 10^{-5} \, {\rm mM pA}^{-1}{\rm ms}^{-1}$} \\
\hline
\end{tabular}
\end{table}

The reduced two-dimensional fast subsystem~\eqref{eq:4fast} is obtained from the five-dimensional fast subsystem~\eqref{eq:7fast} by taking quasi-steady-state assumptions for $m$ and $s$ and setting $h = 1 - 1.08 \, n$.  We use the same parameter values as in Table~\ref{tab:7dpar}, except those specified in Table~\ref{tab:4dpar}, which were adapted to recover the bifurcation structure of the fast subsystem and general slow subsystem behavior as in system~\twoeq{eq:7fast}{eq:7slow}.
\begin{table}[h!]
  \centering
  \caption{\label{tab:4dpar}
    Parameter values used for the functions in the two-dimensional fast subsystem~\eqref{eq:4fast} that are different from those used in the five-dimensional fast subsystem~\eqref{eq:7fast}.}
\begin{tabular}{|r@{ $=$ }r@{.}l|r@{ $=$ }r@{.}l|} 
\hline 
\hline 
$g_{K}$  &   $15$ & $0 \, {\rm nS}$ & $\theta_s$ & $10$ & $0 \, {\rm mV}$ \\
$g_{\rm CAN}$ & $10$ & $0 \, {\rm nS}$ & $k_{\rm CAN}$ & $0$ & $25 \, \mu{\rm M}$ \\
\hline 
\hline 
$\sigma_s$ & $-8$ & $0 \, {\rm mV}$ & $k_{Ca}$ & $60$ & $0 \, {\rm ms}^{-1}$ \\
\multicolumn{3}{|c|}{}   & $k_{\rm IP3}$ & \multicolumn{1}{l}{$1700$} & \multicolumn{1}{r|}{$\mu{\rm M ms}^{-1}$} \\
\hline
\hline
$k$ & $10$ & $0$ & $r_{\rm pump}$ & $1500$ & $0 \, {\rm pA}$ \\
\hline
$\eps$ & $0$ & $005$ & \multicolumn{3}{c|}{} \\
\hline
\end{tabular}
\end{table}
%

\end{document}